\documentclass[11pt]{article}  % Stile del documento
\input{epsf.sty}
\usepackage{amssymb}
\newtheorem{theorem}{Theorem}[section]
\newtheorem{proposition}{Proposition}[section]
\newcommand{\remove}[1]{}
\newenvironment{proof}{
\begin{trivlist}
\item[\hspace{\labelsep}{\bf\noindent Proof. }]
}{$\hfill\Box$\end{trivlist}}
\title{\LARGE\bf A double-ended queue with catastrophes and repairs, and 
a jump-diffusion approximation\footnote{
Paper accepted for publication in {\em Methodology and Computing in Applied Probability}  
\newline
The final publication is available at www.springerlink.com \ (DOI: 10.1007/s11009-011-9214-2)}
}
\author{
Antonio {\bf Di Crescenzo}$^{(1)}$ \\
Virginia {\bf Giorno}$^{(1)}$ \\
Balasubramanian {\bf Krishna Kumar}$^{(2)}$ \\ 
Amelia G. {\bf Nobile}$^{(1)}$\\
$(1)$ Dipartimento di Matematica e Informatica, Universit\`a di Salerno, \\ 
Via Ponte don Melillo, I-84084 Fisciano (SA), Italy\\
email: \{adicrescenzo, giorno, nobile\}@unisa.it\\
$(2)$ Department of Mathematics, Anna University, 
Chennai 600 025, India\\
email: drbkkumar@hotmail.com
}
\date{\empty}
\begin{document}
\maketitle
\begin{abstract}
Consider a system performing a continuous-time random walk on the integers, subject 
to catastrophes occurring at constant rate, and followed by exponentially-distributed 
repair times. After any repair the system starts anew from state zero. 
We study both the transient and steady-state probability laws of the stochastic process 
that describes the state of the system. 
We then derive a heavy-traffic approximation to the model that yields a jump-diffusion 
process. The latter is equivalent to a Wiener process subject to randomly occurring jumps, 
whose probability law is obtained. The goodness of the approximation is finally discussed.


\end{abstract}

\newpage
%   ----------------------------------------------------------------------------
\section{Introduction}
%   ----------------------------------------------------------------------------
The description of numerous types of systems subject to random evolution can 
be rendered in a more effective way by including the eventual occurrence of 
catastrophes. This is well known in the contexts of population dynamics 
(see, for instance, the contributions by Brockwell 1985, Swift 2000 and 2001, and 
Gani and Swift, 2006 and 2007).  
In this research area, Chao and Zheng (2003) considered an immigration 
birth and death process with total catastrophes and studied 
its transient as well as equilibrium behavior. 
Other contributions are due to Economou and Fakinos (2003) and (2008), 
that extend some previous results on birth-death processes with 
catastrophes to the more general cases of continuous-time Markov chains 
with catastrophes resulting from a point process, and of the non-homogeneous 
Poisson process with total or binomial catastrophes. 
Various results on a class of stochastic models for systems subject to random 
regulation via disasters are presented in Stirzaker (2006) and (2007). 
Also the analysis of queueing systems or birth-death processes subject 
to catastrophes has attracted the attention of several investigators (see, 
for instance, Krishna Kumar and Arivudainambi 2000; Krishna Kumar and 
Pavai Madheswari 2003; Di Crescenzo et al.\ 2008). A suitable approach 
adopted by Krinik et al. (2005) to determine the transient probability functions 
of some classical queueing systems with catastrophes uses dual processes, 
randomization and lattice path combinatorics, whereas a case of finite 
birth-death processes with catastrophes was considered in 
Krinik and Mortensen (2007).
\par
Models based on stochastic 
processes in the presence of catastrophes have been recently exploited 
also in another field of mathematical biology, with special reference to the 
description of the interaction between a myosin head and an actin filament. 
Precisely, in Buonocore et al.\ (2009) a continuous-time Markov chain-based model 
and its continuous approximation have been proposed and studied starting from the 
$M/M/1/K$-catastrophe paradigm. 
\par 
Along the research line traced by the above mentioned articles, in this paper we 
consider a different birth-death stochastic model subject to catastrophes. We recall 
that a bilateral birth-death process characterized by constant birth and death rates, 
and subject to jumps to state 0 induced by catastrophes, has been studied recently 
in Di Crescenzo and Nastro (2004) and in Switkes (2004). Here we aim to study an 
extension of such a stochastic model, which performs a continuous-time random 
walk on the integers and, in addition, is subject to the occurrence of 
catastrophes followed by exponentially-distributed repair times. The state-space 
of this model thus consists in the whole set of integers complemented by a 
spurious state, say $F$, which is occupied by the system during the repair times. 
At the end of any repair the system occupies the state zero and starts afresh. 
It should be mentioned that the idea of an $M/M/1$ queueing system subject to 
catastrophes and followed by random repairs has been recently considered in 
Krishna Kumar et al.\ (2007). 
\par
It is worth pointing out that the basic stochastic process of our model, namely the 
bilateral birth-death process with constant birth and death rates, can be viewed for 
instance as describing the state of a double-ended queue in a taxi-passenger system. 
Double-ended queues have been studied e.g.\ in Conolly et al.\ (2002).
Here, denoting by $\{N(t);\;t\geq 0\}$ the state of the considered system, $N(t)=n>0$ 
means that at time $t$ there are $n$ customers in the system, $N(t)=-n<0$ means that 
at time $t$ there are $n$ taxis in the system, $N(t)=0$ means that at time $t$ the 
system is empty, and $N(t)=F$ means that at time $t$ the system is in a repair period. 
\par
The transient analysis of queueing systems is not always easy to be performed, and 
it leads very often to non-manageable formulas. In many cases one is forced to study 
the system under heavy-traffic conditions, in order to obtain at least diffusion 
or continuous approximations. Also systems subject to catastrophes can be approximated 
in such a way, the approximating process being usually a jump-diffusion process. 
This is the case of the contribution given in Di Crescenzo et al.\ (2003), where 
an $M/M/1$ queueing system in the presence of catastrophes is studied and is 
approximated by a Wiener process subject to randomly occurring jumps. A similar 
result will be obtained hereafter for the system under investigation.
\par
This is the plan of the paper. In Section 2 we study the preannounced 
double-ended queue with catastrophes and repairs. By exploting its connection 
with the system in absence of catastrophes we express the state probabilities of 
the queue in terms of integrals of Bessel functions. The mean and the variance are 
also investigated, together with the steady state of the system. Section 3 is devoted 
to develop the approximating procedure that leads to a jump-diffusion process 
of the Wiener type. We thus analyze the approximating process, and obtain its  
transient probability density, mean, variance, and steady-state density. 
In conclusion, in Section 4 we discuss the goodness of the continuous approximation. 
%   ----------------------------------------------------------------------------
\section{Discrete model with catastrophes and repairs}
%   ----------------------------------------------------------------------------
\setcounter{equation}{0}
We consider a system performing a one-dimensional random walk in continuous time 
on the whole set of integers. Let $\{N(t);\;t\geq 0\}$ be the position of the
system at time $t$. The system moves along the real axis such that
$\lambda $ is the rate of moving to the right and $\mu$ is the rate
of moving to the left (see Conolly, 1971). This random walk process
$\{N(t);\;t\geq 0\}$ is also said to be a bilateral birth-death
process, with birth rate $\lambda$ and death rate $\mu$; in any small
interval $(t, t+\Delta t), \,\, \Delta t > 0$, a birth occurs
with probability $ \lambda \Delta t + o(\Delta t)$; a death
occurs with probability $ \mu \Delta t + o(\Delta t)$. It is
clear that in this interval neither a birth nor death takes place
with probability $ 1- (\lambda + \mu ) \Delta t + o(\Delta t)$.
\par
Apart from birth and death processes, catastrophes occur at the
system according to a Poisson process with rate $\nu$, i.e., the
catastrophe occurs at the system in the small interval $(t,t+\Delta t)$ 
with probability $ \nu \Delta t + o(\Delta t)$. Whenever a catastrophe 
occurs at the system, the system goes into
the failure state $F$, i.e., the bilateral birth-death process is
subject to catastrophic failure. The repair times of the failed
system are i.i.d.\ according to an exponential distribution with mean
$\eta^{-1}$. After a repair of the system is completed, the system
immediately returns to state zero and begins to work again 
(``on'' state).
\par
As mentioned in Section 1, $N(t)$ can be seen as the state of a double-ended 
queue subject to disasters and repairs. 
We model the system as a continuous-time Markov process with 
state-space  $S = \{ F\}\cup\mathbb{Z}= \{ F, 0, \pm1, \pm2,...\}$. 
Let $P_n(t) = P[N(t) =n],\,\, n = 0, \pm1, \pm2,...$, be the state probability 
that the system is in state $n$  at time $t$ when the system is working
(``on'' state) and $q(t)$ be the probability that the system at time $t$ is under
repair (``failure'' state).
\par
For the system under investigation, the Chapman-Kolmogorov forward
differential-difference equations for the system state
probabilities $P_n(t) ,\,\, n = 0, \pm1, \pm2,...$, and failure
state probability $q(t)$ can be written as
\begin{eqnarray}
 && \frac{{\rm d} q(t)}{{\rm d}t}=-\eta\,q(t)+\nu\,[1-q(t)],
 \label{equation:systemp1} \\
 && \frac{{\rm d} P_0(t)}{{\rm d}t}=-(\lambda+\mu+\nu)\,P_0(t)+\lambda\,P_{-1}(t)+
 \mu\,P_1(t) +\eta\,q(t),
 \\
 && \frac{{\rm d} P_n(t)}{{\rm d}t}=-(\lambda+\mu+\nu)\,P_n(t)+\lambda\,P_{n-1}(t)+\mu\,P_{n+1}(t),
 \qquad n = \pm 1,\pm 2, \pm 3, ...\,\,.
 \label{equation:systemp4}
\end{eqnarray}
Without loss of generality, assume that initial position of the
system at time  $t = 0$ is zero, i.e., 
\begin{equation}
  q(0)=0 \quad \hbox{and}  \quad P_n(0)=\cases{
 1 & if $n=0$,\cr
 0 & otherwise.}
 \label{equation:initial}
\end{equation}
%
%   ----------------------------------------------------------------------------
\subsection{Transient distribution and relation with the process in absence of catastrophes}
%   ----------------------------------------------------------------------------
Let us now introduce the so-called randomized random walk $\{\widetilde N(t);\;t\geq 0\}$ 
with parameters $\lambda$ and  $\mu$. It is a bilateral birth-death 
process over $\mathbb{Z}$ characterized by birth rate $\lambda$ and death rate 
$\mu$, and describes the state of the system under consideration in the 
absence of catastrophes. We recall that the state probabilities of process 
$\widetilde N(t)$ conditional on $\widetilde N(0)=0$ are given by 
(see, for instance, Baccelli and Massey 1989 or Conolly 1971) 
\begin{equation}
 \widetilde P_n(t)= {\rm e}^{-(\lambda+\mu) t}\beta^{n}
 \,I_{n}\left(\alpha t\right),
 \qquad t\geq 0, \;\;  n\in\mathbb{Z},
 \label{equation:Skellam}
\end{equation}
where 
\begin{equation}
 \alpha = 2 \sqrt{\lambda \mu}, 
 \qquad 
 \beta = \sqrt{\frac{\lambda}{\mu}} 
\label{equation:alfabeta}
\end{equation}
and $I_n(\cdot)$ is the modified Bessel function of the first kind of order $n$. 
We remark that the distribution given in (\ref{equation:Skellam}) is often named 
Skellam distribution. 
\par
Hereafter we express the state probabilities $P_n(t)$ in terms of suitable 
quantities concerning the bilateral birth-death process in absence of catastrophes 
$\widetilde N(t)$. Indeed, from probabilistic arguments it is not hard to see that the 
probabilities $P_n(t)$ can be expressed in terms of (\ref{equation:Skellam}) as follows: 
\begin{equation}
 P_n(t)= e^{-\nu t} \,\widetilde P_n(t)
 +\eta
 \int_0^t q(\tau)\, e^{-\nu\,(t-\tau)}
 \,\widetilde P_n(t-\tau)\,{\rm d}\tau,
 \qquad t\geq 0, \;\;  n\in\mathbb{Z}.
 \label{equation:PedR}
\end{equation}
The first term in the right-hand-side of (\ref{equation:PedR}) 
is the probability that the system goes from 0 to $n$ without occurrence 
of catastrophes in $(0,t)$, whereas the second term represents the probability 
that at time $\tau\in(0,t)$ the system enters state $0$ due to a repair 
and then goes from 0 to $n$ without occurrence of catastrophes in $(\tau,t)$. 
\begin{theorem} \label{theorem:prob}
For all $t\geq 0$ we have
\begin{equation}
q(t) = \frac{\nu}{\eta + \nu} \left[1 - e^{-(\eta + \nu)t}\right], 
\label{equation:q(t)}
\end{equation}
%
% \begin{eqnarray}
% &&
% P_0(t) = e^{-(\lambda +\mu +\nu)t} I_0(\alpha t) 
% + \frac{\eta\nu}{\eta +\nu}\int_0^t e^{-(\lambda +\mu +\nu)u} I_0(\alpha u) {\rm d}u
% \nonumber \\
% && \hspace{1.5cm}- \frac{\eta \nu}{\eta +\nu}e^{-\nu t}\int_0^t
% e^{-\eta u} e^{-(\lambda +\mu )(t-u)} I_0(\alpha (t-u)) {\rm d}u, 
% \label{equation:p0(t)}
% \end{eqnarray}
%
\begin{eqnarray}
&&
P_n(t) = \beta^n \bigg\{ e^{-(\lambda +\mu +\nu)t} I_n(\alpha t) 
+\frac{\eta \nu}{\eta+\nu}\bigg[
\int_0^t e^{-(\lambda +\mu +\nu) u} I_n(\alpha u) {\rm d}u 
\nonumber   \\
&& \hspace{1.5cm}
-e^{-\nu t} \int_0^t e^{-\eta \tau - (\lambda +\mu)(t-\tau)} I_n(\alpha (t-\tau)) {\rm d}\tau\bigg]\bigg\},
\qquad n \in \mathbb{Z},
\label{equation:pn(t)}
\end{eqnarray}
where $\alpha$ and $\beta$ are defined in (\ref{equation:alfabeta}). 
\end{theorem}
\begin{proof}
\par
Making use of Eqs.\  (\ref{equation:systemp1}) and recalling the first of 
(\ref{equation:initial}) we obtain (\ref{equation:q(t)}). 
Furthermore, substituting (\ref{equation:Skellam}) and (\ref{equation:q(t)}) in 
Eq.\ (\ref{equation:PedR}) the transient distribution (\ref{equation:pn(t)}) follows. 
\end{proof}
\par
Some plots of probabilities (\ref{equation:pn(t)}) are shown in Figure 1. 
The expressions given in Theorem \ref{theorem:prob} completely determine 
all state transient probabilities of the system and the failure state probability 
$q(t)$ for the bilateral birth-death process $N(t)$ with catastrophes and repairs. 
%
%--------------------------------------------------------------------
\begin{figure}[b]   %%%%%%%%%%% FIGURA  probabilitˆ
%\vspace{20mm}
\begin{center}
\epsfxsize=7.3cm
\epsfbox{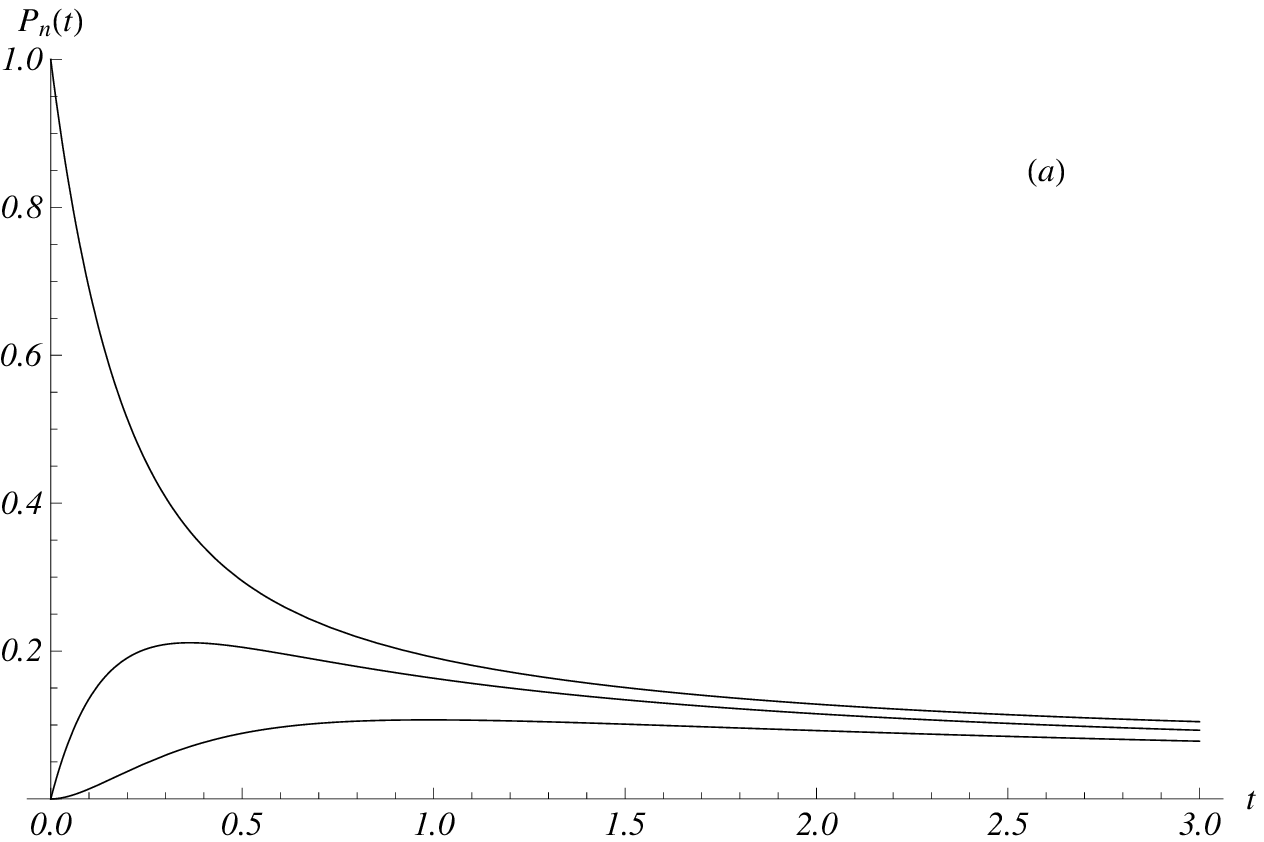} 
$\;\;$
\epsfxsize=7.3cm
\epsfbox{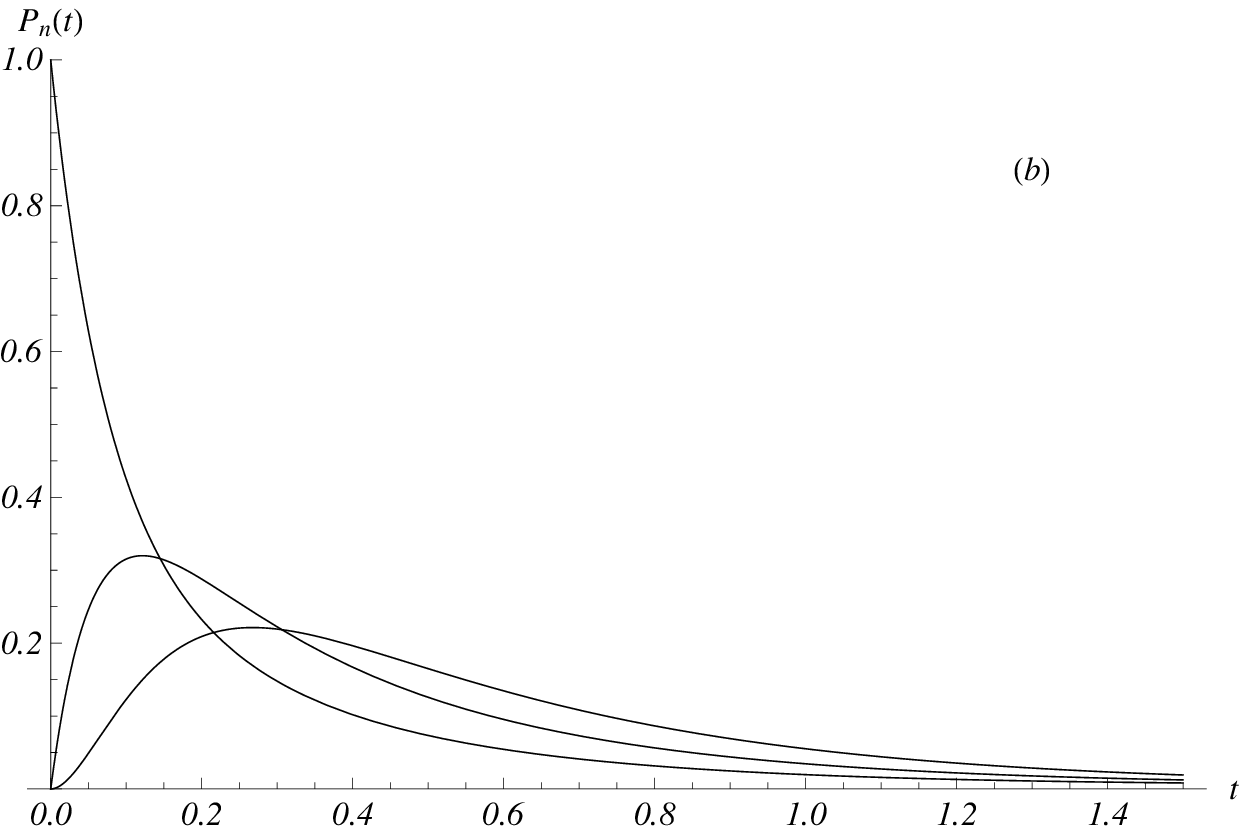}
\vspace{-0.2cm}
\caption{Plots of $P_n(t)$ for $n=0,1,2$ (from top to bottom near the origin), 
where $\mu = 2$, $\nu = 0.1$, $\eta = 1$, and (a) $\lambda = 2$, (b) $\lambda = 8$.}
\end{center}
\end{figure}
%--------------------------------------------------------------------
%
\par
Let us now denote by $\widetilde T_{0,n}$ the first-passage 
time from state $0$ to state $n$ $( n = \pm 1, \pm 2, \pm 3, \ldots)$ of process 
$\widetilde N(t)$, and by $\widetilde g_{0,n}(t)$ the corresponding probability 
density function. It is well known that (see Abate et al.\ 1991, for instance), for all $t>0$  
\begin{equation}
 \widetilde g_{0,n} (t) 
 = \frac{|n|}{t} \widetilde P_n(t)
 = \frac{|n|}{t} e^{-(\lambda +\mu) t} \beta^n
 I_n(\alpha t),\qquad n = \pm 1, \pm 2, \pm 3, \ldots\,\,.
 \label{equation:fptd}
\end{equation}
\par 
Due to the Markov property, the following renewal-type equation holds
for the system under consideration, for $t>0$:
\begin{equation}
 P_n(t) = \int_0^t P_0(u) e^{-\nu (t-u)}\widetilde g_{0,n}(t-u)\, {\rm d}u,
 \qquad n = \pm 1, \pm 2, \pm 3, \ldots \,\,.
 \label{equation:integrpn(t)}
\end{equation}
Hence, Eq.\ (\ref{equation:integrpn(t)}) expresses 
that any sample-path of $N(t)$ going from 0 to state $n$, at an eventual time 
$u\in(0,t)$ visits the state 0 and then goes from 0 to $n$ without occurrence 
of catastrophes in $(u,t)$. Hereafter we provide an alternative 
expression of the transition probabilities. 
\begin{proposition}
For all $t\geq 0$ the following relation holds:
$$
P_n(t) = |n| \beta^n \int_0^t e^{-(\lambda +\mu +\nu) (t-u)}
\frac{I_n(\alpha (t-u))}{t-u} P_0(u) {\rm d}u, 
\qquad n = \pm 1, \pm 2, \pm 3, ...
$$
\end{proposition}
\begin{proof}
It immediately follows by making use of  (\ref{equation:fptd}) in 
(\ref{equation:integrpn(t)}). 
\end{proof}
\par
In the following, for any function $f(t)$, $t\geq 0$, we denote by 
$f^*(z)=\int_0^{\infty} e^{-z t}\,f(t)\,{\rm d}t$, $z>0$, its Laplace transform. 
\par
By adopting customary methodologies, such as those described in   
Jain et al.\ (2007) or Parthasarathy and Lenin  (2004), we can evaluate the 
Laplace transform of the transition probabilities $P_n(t)$. Indeed, by taking 
Laplace transforms in the system of equations 
(\ref{equation:systemp1})-(\ref{equation:systemp4}), after some calculations we obtain:
\begin{equation}
P_0^*(z) = \frac{1}{\sqrt{(z+\lambda +\mu +\nu)^2-4 \lambda
\mu}}\,+\,\frac{\eta \nu}{z(z+\eta + \nu) \sqrt{(z+\lambda +\mu
+\nu)^2-4 \lambda \mu}}
\label{equation:p*p0z}
\end{equation}
and
\begin{equation}
P_n^*(z) = \left\{
\begin{array}{ll}
P_0^*(z)\,[\psi_2(z)]^n, & n=1,2,3,... \\
P_0^*(z)\,[\psi_1(z)]^n, & n=-1,-2,-3,...  
\end{array}
\right.
\label{equation:p*pnz}
\end{equation}
where 
$$
 \psi_1(z), \psi_2(z) = \frac{z+\lambda +\mu +\nu \pm
\sqrt{(z+\lambda +\mu +\nu)^2-4 \lambda \mu}}{2\mu},
\qquad \psi_1(z)> \psi_2(z)
$$ 
are the roots of quadratic equation 
$\mu \chi^2(z) -(z+\lambda +\mu + \nu)  \chi(z) + \lambda = 0$. 
Expressions (\ref{equation:p*p0z}) and (\ref{equation:p*pnz}) 
will be used in Section 4 in order to discuss the goodness of the 
approximation of $N(t)$ by a suitable jump-diffusion process $X(t)$. 
%   ----------------------------------------------------------------------------
\subsection{Steady-state distribution}
%   ----------------------------------------------------------------------------
In this section, we shall investigate the structure of the steady-state 
probabilities and the failure state probability of process $N(t)$:
$$
 \pi_n:=\lim_{t\to \infty} P_n(t), \qquad
 q:=\lim_{t\to \infty} q(t).
$$
\begin{theorem}
The steady-state probabilities of the system state $\left\{\pi_n: n
= 0, \pm 1, \pm 2, \ldots\right\} $ and the failure probability $q$
of the bilateral birth-death process with catastrophes are obtained
as, for $\nu > 0$ and $\eta >0$,
\begin{eqnarray}
 &&
q = \frac{\nu}{\eta +\nu},
 \label{equation:qstst}
\\
&&
 \pi_0
 = (1-q) \frac{\nu}{\sqrt{(\lambda+\mu+\nu)^2-4\lambda\mu}},
 \label{equation:pi0}
\\
 && \pi_n=\left[\frac{(\lambda+\mu+\nu)
 -\sqrt{(\lambda+\mu+\nu)^2-4\lambda\mu}}{ 2\mu}\right]^n
 \pi_0,
 \qquad n=1,2,\ldots
 \label{equation:pin}
\end{eqnarray}
and
\begin{eqnarray}
 \pi_{-n}=\left[\frac{(\lambda+\mu+\nu)
 -\sqrt{(\lambda+\mu+\nu)^2-4\lambda\mu}}{ 2\lambda}\right]^n
 \pi_0  \qquad n=1,2,\ldots.
 \label{equation:aspimn}
\end{eqnarray}
\end{theorem}
\begin{proof}
The steady-state probabilities of the system can be obtained by letting $t\to +\infty$ in 
the system of equations (\ref{equation:systemp1})-(\ref{equation:systemp4}). The 
resulting balance equations yield a system of recurrence relations of degree 2 with 
known coefficients, whose solution can be obtained by standard methods. 
Eqs.\ (\ref{equation:qstst})-(\ref{equation:aspimn}) then follow. 
\end{proof}
%
%   ----------------------------------------------------------------------------
\subsection{Moments}
%   ----------------------------------------------------------------------------
Let us now focus on the mean of the bilateral birth-death process 
with catastrophes and repairs, and set 
$$
 m_N(t) = E[N(t)\cdot {\bf 1}_{\{N(t)\neq F\}}\,|\,N(0)=0], \qquad t\geq 0.
$$ 
Then, from (\ref{equation:PedR}) we have 
$$
 m_N(t)= e^{-\nu t} \,\widetilde m(t)
 +\eta
 \int_0^t q(\tau)\, e^{-\nu\,(t-\tau)}
 \, m_{\widetilde N}(t-\tau)\,{\rm d}\tau,
 \qquad t\geq 0,  
$$
where $m_{\widetilde N}(t):=E[\widetilde N(t)\,|\,\widetilde N(0)=0]=(\lambda-\mu)t$ 
is the mean of the randomized random walk with initial state 0. 
Hence, recalling (\ref{equation:q(t)}) after some calculations we obtain 
\begin{eqnarray}
m_N(t) = \frac{(\lambda -\mu)\eta}{(\eta+\nu) \nu}
  \left\{1 - e^{-\nu t} 
 +\frac{\nu^2}{\eta^2} e^{-\nu t}(1 -   e^{-\eta  t}) \right\}.
\label{eq:medfinale}
\end{eqnarray}
Due to (\ref{eq:medfinale}), if  $\lambda >\mu$ the following holds: 
\par
-- if $\eta\geq \nu$ then $m_N(t)$ is increasing for all $t>0$, 
\par
-- if $\eta< \nu$ then $m_N(t)$ is increasing for $t<\bar t$ and is decreasing 
for $t>\bar t$, where 
$$
 \bar t :={1\over \eta}\log{\nu \over \nu-\eta}.
$$ 
Clearly, if $\lambda <\mu$ the monotonicity of $m_N(t)$ is reversed. 
Some plots of $m_N(t)$ are shown in Figure 2. 
\par
Let us now introduce the variance of the bilateral birth-death process 
with catastrophes and repairs:
$$
 V_N(t) = Var[N(t)\cdot {\bf 1}_{\{N(t)\neq F\}}\,|\,N(0)=0], \qquad t\geq 0.
$$  
By adopting a similar procedure, making use of (\ref{equation:PedR}), 
(\ref{eq:medfinale}) and recalling that 
$Var[\widetilde N(t)\,|\,\widetilde N(0)=0]=(\lambda+\mu)t$, 
it is not hard to obtain $V_{N}(t)$. Some plots of the variance of $N(t)$
are given in Figure 3. 
\par
We now discuss the asymptotic mean and variance of process $N(t)$. 
From (\ref{eq:medfinale}) it is not hard to see that,  for $v>0$, 
\begin{eqnarray}
 \lim_{t\to +\infty}m_N(t) 
 =   \frac{(\lambda -\mu) \eta}{(\eta+\nu)\nu} 
 =   (1-q)\,\frac{\lambda -\mu}{ \nu }.
 \label{eq:medasfinale}
\end{eqnarray}
We point out that the presence of catastrophes and repairs may be seen as a 
regulatory effect to the system. Indeed, if $\nu=0$ the system does not admit a 
steady-state distribution, whereas if $\nu>0$ then the right-hand-side of (\ref{eq:medasfinale}) is finite. 
Finally, straightforward calculations lead to the following asymptotic variance,  for $v>0$: 
\begin{equation}
 \lim_{t\to +\infty}V_N(t) 
 ={(\lambda+\mu)\eta\over (\eta+\nu)\nu}
 +{(\lambda-\mu)^2 \eta\over (\eta+\nu)^2\nu^2}(2\nu+\eta).
 \label{equation:varN}
\end{equation}
%
%--------------------------------------------------------------------
\begin{figure}[b]   %%%%%%%%%%% FIGURA  media
%\vspace{20mm}
\begin{center}
\epsfxsize=7.3cm
\epsfbox{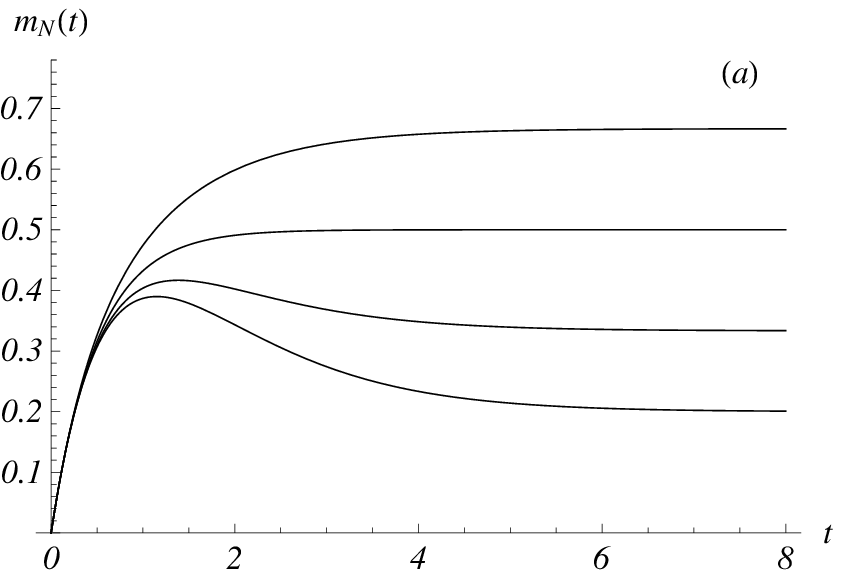} 
$\;\;$
\epsfxsize=7.3cm
\epsfbox{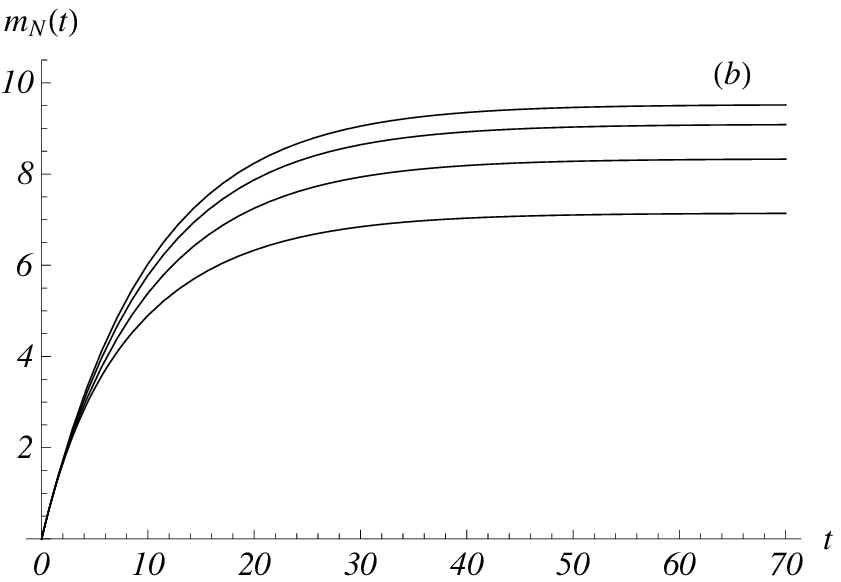}
\vspace{-0.2cm}
\caption{Plots of $m_N(t)$ for $\lambda-\mu=1$, 
 $\eta=0.25$, $0.5$, $1$, $2$ (from bottom to top), and (a) $\nu=1$, (b)  $\nu=0.1$.}
\end{center}
\end{figure}
%--------------------------------------------------------------------
%
%--------------------------------------------------------------------
\begin{figure}[t]   %%%%%%%%%%% FIGURA  varianza
%\vspace{20mm}
\begin{center}
\epsfxsize=7.3cm
\epsfbox{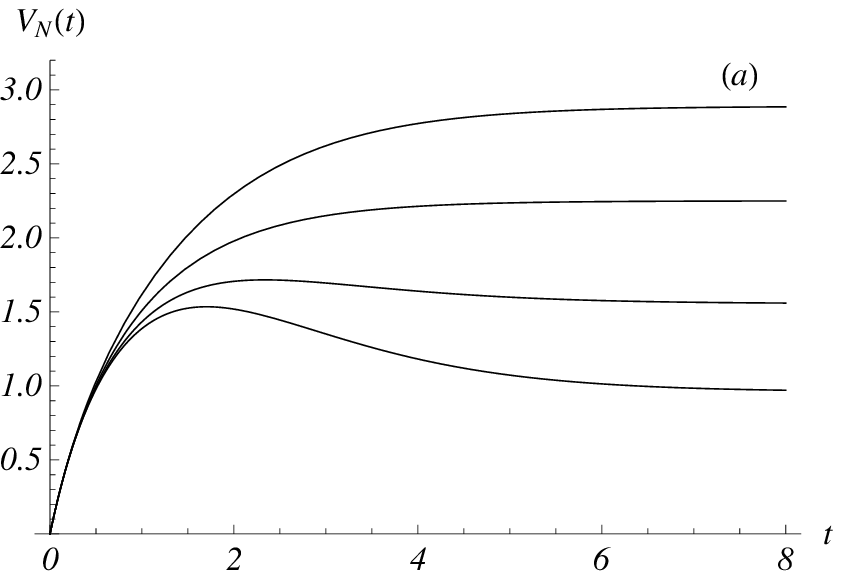} 
$\;\;$
\epsfxsize=7.3cm
\epsfbox{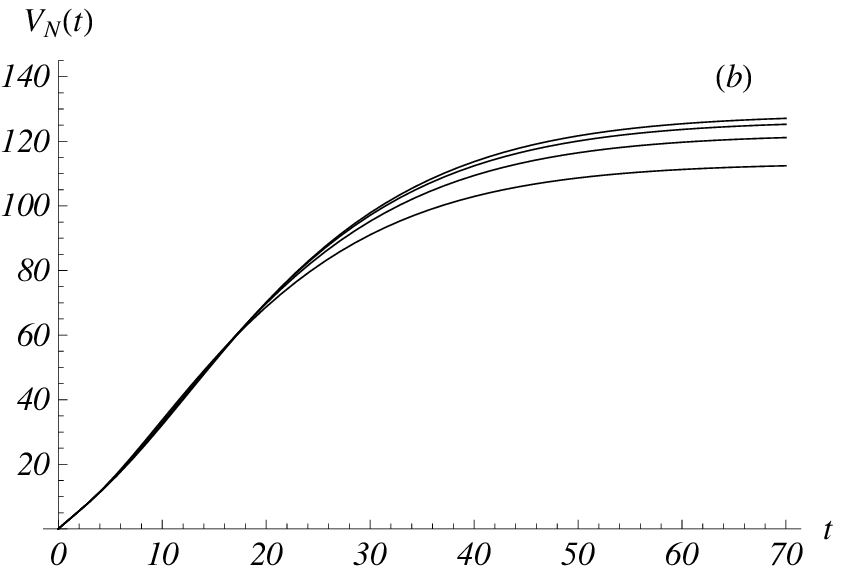}
\vspace{-0.2cm}
\caption{Plots of $V_N(t)$ for the same parameters of Figure 2, with $\lambda=2$.}
\end{center}
\end{figure}
%--------------------------------------------------------------------
%
%   ----------------------------------------------------------------------------
\section{Jump-diffusion approximation}
%   ----------------------------------------------------------------------------
\setcounter{equation}{0}
Let us now construct a continuous approximation starting from the process 
$N(t)$. We shall make use of a scaling similar to one succesfully used elsewhere 
in queueing theory contexts (see Giorno et al. 1986, 1987; Di Crescenzo and
Nobile 1995). We recall that a popular paper dealing with Brownian motion with 
jumps that approximates a standard single-server queue with random vacations 
is due to Kella and Whitt (1990), whereas Kimura (2004) proposed an useful 
review of some basic concepts and issues in diffusion modeling, and a 
bibliographical guide to diffusion models for queues that are typically found 
in computer/communication systems. 
\par
Let us consider the Markov process $\{N^a_{\varepsilon}(t);\;t\geq 0\}$, 
having state space
$\{F\}\cup \{\ldots,-\varepsilon,0,\varepsilon,\ldots\}$, with
$\varepsilon>0$. For all $t>0$, the transitions of
$N^a_{\varepsilon}(t)$ are governed by the following rates:
\begin{eqnarray*}
 && \lim_{\Delta t\downarrow 0}{1\over \Delta t}
 P\{N^a_{\varepsilon}(t+\Delta t)=(n+1)\varepsilon\,|\,N^a_{\varepsilon}(t)=n\varepsilon\}
 = {\widehat \lambda\over \varepsilon}+{\sigma^2\over 2\varepsilon^2},
 \qquad n\in{\mathbb Z},
 \\
 && \lim_{\Delta t\downarrow 0}{1\over \Delta t}
 P\{N^a_{\varepsilon}(t+\Delta t)=(n-1)\varepsilon\,|\,N^a_{\varepsilon}(t)=n\varepsilon\}
 = {\widehat \mu\over \varepsilon}+{\sigma^2\over 2\varepsilon^2},
 \qquad n\in{\mathbb Z},
 \end{eqnarray*}
where $\widehat \lambda$, $\widehat \mu$ and $\sigma$ are positive
parameters. Moreover, the process $\{N^a_{\varepsilon}(t);\;t\geq 0\}$ 
is still subject to catastrophes arriving with rate $\nu$.
Repairs are still exponentially distributed with rate $\eta$, and
the system enters state 0 just after a repair. Hence, the repair and
failure time distributions are not affected by the scaling procedure
introduced above, so that
\begin{eqnarray*}
 && \lim_{\Delta t\downarrow 0}{1\over \Delta t}
 P\{ N^a_{\varepsilon}(t+\Delta t)= F \,|\,N^a_{\varepsilon}(t)=n\varepsilon\}
 =\nu,
 \qquad n\in{\mathbb Z},
 \\
 && \lim_{\Delta t\downarrow 0}{1\over \Delta t}
 P\{N^a_{\varepsilon}(t+\Delta t)=0\,|\,N^a_{\varepsilon}(t) = F\}
 = \eta,
 \qquad n\in{\mathbb Z},
\end{eqnarray*}
where $\eta$ and $\nu$ are the same parameters of $N(t)$.
Essentially, the integers of the state-space of $N(t)$ have been rescaled by
a factor $\varepsilon$, and the following substitutions  have been performed:
\begin{equation}
 \lambda={\widehat \lambda\over \varepsilon}+{\sigma^2\over 2\varepsilon^2},
 \qquad
 \mu={\widehat \mu\over \varepsilon}+{\sigma^2\over 2\varepsilon^2}.
 \label{equation:parameters}
\end{equation}
For $t\geq 0$ the transient behaviour of $N^a_{\varepsilon}(t)$ is described by
the following probabilities:
\begin{equation}
 P^a_n(t)
 =P\{N^a_{\varepsilon}(t)=n\,\varepsilon\},
 \qquad n\in \mathbb Z.
 \label{equation:probeps}
\end{equation}
As before, assuming that $P\{N^a_{\varepsilon}(0)=0\}=1$, we have $P^a_0(0)=1$. 
Moreover, we note that
$$
 \sum_{n\in\mathbb Z} P^a_n(t)=1-q(t),
$$
where $q(t)$ is the probability that $N^a_{\varepsilon}(t)$ is in state $F$, and is 
given in  (\ref{equation:q(t)}).
\par
Let us now consider a continuous-time Markov process $\{X(t);\;t\geq 0\}$ 
having state space $\{F\}\cup(-\infty,\infty)$. This is a jump-diffusion process that  
describes the motion of a particle that starts from the origin at time $0$, then 
behaves as a Wiener process with drift $\hat \lambda - \hat \mu$ and infinitesimal 
variance $\sigma^2$, until it jumps to state $F$ with rate $\nu$. Subsequently it jumps 
again to the initial state 0, with rate $\eta$, and diffuses afresh. The probability law 
of $X(t)$ possesses an absolutely continuous component 
$$
 f(x,t\,|\,0)=\lim_{h\downarrow 0}{1\over h}
 P\{x\leq X(t)<x+h\,|\,X(0)=0\}, \qquad x\in\mathbb{R}, \;\; t\geq 0
$$ 
and a discrete component $P\{X(t)=F\,|\,X(0)=0\}$. The latter is easily seen to be 
identical to the probability $q(t)$, obtained in (\ref{equation:q(t)}). 
\par
The continuous approximation relies on the convergence of $N^a_{\varepsilon}(t)$ 
to process $X(t)$ as $\varepsilon$ tends to 0. Probabilities 
(\ref{equation:probeps}) are near to $f(n\,\varepsilon,t\,|\,0)\,\varepsilon$ 
when $\varepsilon$ is close to zero, and the correspondence between such 
quantities becomes increasingly better as $\varepsilon$ approaches zero.
We thus follow a customary procedure based on the substitution of $P^a_n(t)$ with 
$f(n\,\varepsilon,t\,|\,0)\,\varepsilon$ in the Chapman-Kolmogorov forward 
differential-difference equations of $N^a_{\varepsilon}(t)$. After setting $x=n\,\varepsilon$ 
and expanding $f$ as Taylor series, we pass to the limit as $\varepsilon\downarrow 0$. 
By taking positions (\ref{equation:parameters}) into account and assuming that 
$\Delta t$ is proportional to $\varepsilon^2$, after some calculations we obtain:
\begin{equation}
 {\partial\over\partial t} f(x,t\,|\,0)
 =-\nu\,f(x,t\,|\,0)
 -(\widehat\lambda-\widehat\mu)\,{\partial\over\partial x} f(x,t\,|\,0)
 +{\sigma^2\over 2}\,{\partial^2\over\partial x^2} f(x,t\,|\,0),
 \qquad x\neq 0,
 \label{equation:pdiff}
\end{equation}
to be solved with condition
\begin{equation}
 \int_{-\infty}^{\infty}f(x,t\,|\,0)\,{\rm d}x=1-q(t),
 \qquad t\geq 0.
 \label{equation:cont}
\end{equation}
Equation (\ref{equation:pdiff}) is the Fokker-Planck equation for a
Wiener process  perturbed by a leaky term $-\nu\,f(x,t\,|\,0)$.
Note that the initial condition $P^a_0(0)=1$ now becomes
\begin{equation}
 \lim_{t\downarrow 0} f(x,t\,|\,0)=\delta(x),
 \label{equation:init}
\end{equation}
where $\delta(x)$ is the Dirac delta-function.
\par
The above sketched procedure indicates that the  process $N^a_{\varepsilon}(t)$
converges in some sense to the process $X(t)$ as $\varepsilon$ tends to 0. Hence, 
statistical indexes of $X(t)$ are suitable to determine quantities useful to describe $N(t)$ 
when $\varepsilon$ is very small. Indeed, under position  (\ref{equation:parameters}) we have
$$
 P\{N(t)<n\,|\,N(0)=0\}\simeq
 P\{X(t)<n\,\varepsilon\,|\,X(0)=0\};
$$
recalling that $x=n\,\varepsilon$, this approximation is expected
to improve as $\varepsilon$ goes to zero and as $|n|$ grows larger.
%   ----------------------------------------------------------------------------
\subsection{Transient and steady-state distributions}
%   ----------------------------------------------------------------------------
In order to analyse the  process $X(t)$ hereafter we obtain the function 
$f^*(x,z\,|\,0)$, that is the Laplace transform of density $f(x,t\,|\,0)$. 
Recalling (\ref{equation:init}), from (\ref{equation:pdiff}) we then have
\begin{equation}
 (z+\nu)\,f^*(x,z\,|\,0)
 =-(\widehat\lambda-\widehat\mu)\,{\partial\over\partial x}f^*(x,z\,|\,0)
 +{\sigma^2\over 2}\,{\partial^2\over\partial x^2}f^*(x,z\,|\,0),
 \qquad x\neq 0,
 \label{equation:eqphi}
\end{equation}
whereas condition (\ref{equation:cont}) becomes
\begin{eqnarray}
 \int_{-\infty}^{\infty}f^*(x,z\,|\,0)\,{\rm d}x={z+\eta\over z\,(z+\eta+\nu)}.
 \label{equation:intphi}
\end{eqnarray}
\begin{theorem}
For $z>0$ we have
\begin{eqnarray}
 f^*(x,z\,|\,0)
 \!\!\! &=& \!\!\! {1\over \sqrt{(\widehat\lambda-\widehat\mu)^2 +2\sigma^2(z+\nu)}}\,
 {(z+\nu)(z+\eta)\over z\,(z+\eta+\nu)}
 \nonumber \\
 &\times & \!\!\!
 \exp\left\{{\widehat\lambda-\widehat\mu\over\sigma^2}\,x
 -{\sqrt{(\widehat\lambda-\widehat\mu)^2 +2\sigma^2(z+\nu)}\over\sigma^2}\,|x|\right\}
 \qquad (x\in\mathbb R).
 \label{equation:solphi}
\end{eqnarray}
\end{theorem}
\begin{proof}
By continuity of $f^*(x,z\,|\,0)$ in $x=0$
the solution of (\ref{equation:eqphi}) takes the form
\begin{equation}
 f^*(x,z\,|\,0)=\cases{
 A\, e^{w_2 \,x} & if $x>0$ \cr
 A\, e^{w_1 \,x} & if $x<0$
}
 \label{equation:psolphi}
\end{equation}
where
\begin{equation}
 w_1,w_2={\widehat\lambda-\widehat\mu\,\pm\,
 \sqrt{(\widehat\lambda-\widehat\mu)^2 +2\sigma^2(z+\nu)}\over\sigma^2}
 \qquad (w_2<0<w_1)
 \label{equation:wespr}
\end{equation}
are solutions of equation
$$
 \sigma^2 \,w^2-2(\widehat\lambda-\widehat\mu)\,w-2(z+\nu)=0.
$$
Making use of (\ref{equation:psolphi}) in (\ref{equation:intphi}) we thus obtain
\begin{equation}
 A={1\over \sqrt{(\widehat\lambda-\widehat\mu)^2 +2\sigma^2(z+\nu)}}\,
 {(z+\nu)(z+\eta)\over z\,(z+\eta+\nu)}.
 \label{equation:Aespr}
\end{equation}
Hence, Eq.\ (\ref{equation:solphi}) finally follows from
(\ref{equation:psolphi}), (\ref{equation:wespr}) and (\ref{equation:Aespr}).
\end{proof}
\par 
We are now able to obtain an expression of $f(x,t\,|\,0)$, that is the analogous 
of (\ref{equation:PedR}) for the discrete system.
\begin{theorem}
For all $x\in\mathbb R$ and $t\geq 0$ we have
\begin{equation}
 f(x,t\,|\,0)= e^{-\nu t} \,\widetilde f(x,t\,|\,0)
 +\eta
 \int_0^t q(\tau)\, e^{-\nu\,(t-\tau)}
 \,\widetilde f(x,t-\tau\,|\,0)\,{\rm d}\tau,
 \label{equation:f(x,t)}
\end{equation}
where $\widetilde f(x,t\,|\,0)$ is the transition density of a Wiener process $\{\widetilde X(t);\;t\geq 0\}$
having drift $\widehat\lambda-\widehat\mu$ and infinitesimal variance $\sigma^2$.
\end{theorem}
\begin{proof}
First of all we note that Eq.\ (\ref{equation:solphi}) can be expressed as
\begin{eqnarray}
 f^*(x,z\,|\,0)
 \!\!\! &=& \!\!\!
  \left(1+{\eta\,\nu\over \eta+\nu}\,{1\over z}
 -{\eta\,\nu\over \eta+\nu}\,{1\over z+\eta+\nu}\right)
 \exp\left\{{\widehat\lambda-\widehat\mu\over\sigma^2}\,x\right\}
 \nonumber \\
 &\times & \!\!\!
 {1\over \sqrt{(\widehat\lambda-\widehat\mu)^2 +2\sigma^2(z+\nu)}}
 \exp\left\{-{\sqrt{(\widehat\lambda-\widehat\mu)^2 +2\sigma^2(z+\nu)}\over\sigma^2}\,|x|\right\}
 \quad (x\in\mathbb R).\qquad
 \label{equation:phialtern}
\end{eqnarray}
Let us now consider a Wiener process having drift $\sqrt{(\widehat\lambda-\widehat\mu)^2+2\sigma^2\nu}$
and infinitesimal variance $\sigma^2$; its transition density $f_{\scriptscriptstyle W}(x,t\,|\,0)$ has
Laplace transform
$$
 f_{\scriptscriptstyle W}^*(x,z\,|\,0)
 ={1\over \sqrt{(\widehat\lambda-\widehat\mu)^2 +2\sigma^2(z+\nu)}}\,
 \exp\left\{{\sqrt{(\widehat\lambda-\widehat\mu)^2 +2\sigma^2 \nu}\over\sigma^2}\,x
 -{\sqrt{(\widehat\lambda-\widehat\mu)^2 +2\sigma^2(z+\nu)}\over\sigma^2}\,|x|\right\},
$$
so that
\begin{equation}
 {\exp\left\{-\displaystyle{|x|\over\sigma^2}\,\sqrt{(\widehat\lambda-\widehat\mu)^2 +2\sigma^2(z+\nu)}\right\}
 \over \sqrt{(\widehat\lambda-\widehat\mu)^2 +2\sigma^2(z+\nu)}}
 = \exp\left\{-{\sqrt{(\widehat\lambda-\widehat\mu)^2 +2\sigma^2 \nu}\over\sigma^2}\,x\right\}
 f_{\scriptscriptstyle W}^*(x,z\,|\,0)
 \label{equation:relfstar}
\end{equation}
Making use of  Eq.\ (\ref{equation:relfstar}), thus Eq.\ (\ref{equation:phialtern}) yields
\begin{equation}
f^*(x,z\,|\,0)
 =e^{z_2\,x}\left(1+{\eta\,\nu\over \eta+\nu}\,{1\over z}
 -{\eta\,\nu\over \eta+\nu}\,{1\over z+\eta+\nu}\right)
 f_{\scriptscriptstyle W}^*(x,z\,|\,0)
 \qquad (x\in\mathbb R),
 \label{equation:phiugftilde}
\end{equation}
where we have set
$$
 z_2={\widehat\lambda-\widehat\mu-\sqrt{(\widehat\lambda-\widehat\mu)^2 +2\sigma^2 \nu}
 \over \sigma^2}\,\,.
$$
Hence, by taking the inverse Laplace transform in (\ref{equation:phiugftilde}), we come to
$$
 f(x,t\,|\,0)
 =e^{z_2\,x}
 \left\{f_{\scriptscriptstyle W}(x,t\,|\,0)
 +{\eta\,\nu\over \eta+\nu}\int_0^t \left[1-e^{-(\eta+\nu)(t-\tau)}\right]
 \,f_{\scriptscriptstyle W}(x,\tau\,|\,0)\,{\rm d}\tau\right\}.
$$
Finally, recalling (\ref{equation:q(t)}) and making use of  identity
$$
 e^{z_2\,x}\,f_{\scriptscriptstyle W}(x,t\,|\,0)
 = e^{-\nu t}\,\widetilde f(x,t\,|\,0),
$$
Eq.\  (\ref{equation:f(x,t)}) follows.
\end{proof}
\par
Some plots of density $f(x,t\,|\,0)$ obtained from (\ref{equation:f(x,t)}) are 
given in Figure 4. 
%--------------------------------------------------------------------
\begin{figure}[b]   %%%%%%%%%%% FIGURA  f(x,t)
%\vspace{20mm}
\begin{center}
\epsfxsize=7.3cm
\epsfbox{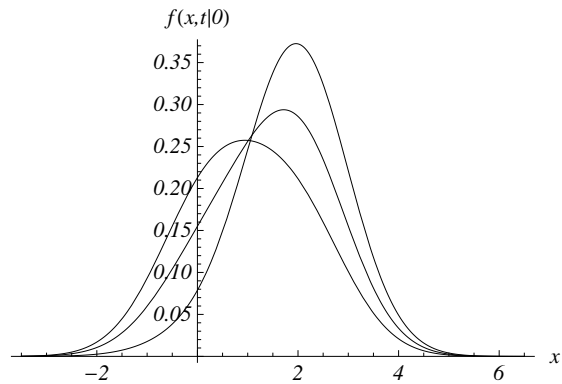}
\vspace{-0.2cm}
\caption{Plots of density (\ref{equation:f(x,t)}),  for $\widehat\lambda-\widehat\mu=2$, 
 $\eta=1$, $\sigma^2=1$, $t=1$, and  $\nu=1$, $0.5$, $0.1$ (from left to right). The probability 
 mass is respectively $0.5677$, $0.7410$, $0.9394$,  due to (\ref{equation:cont}).}
\end{center}
\end{figure}
%--------------------------------------------------------------------
\par 
Hereafter we obtain another expression of $f(x,t\,|\,0)$ which is the analogous of 
(\ref{equation:integrpn(t)}) under the continuous approximation. It involves the 
first-passage-time density of the Wiener process $\widetilde X(t)$ with drift 
$\widehat\lambda-\widehat\mu$ and infinitesimal variance $\sigma^2$, given by: 
\begin{equation}
 \widetilde g(x,t\,|\,x_0) 
 = {|x-x_0|\over t}\,\widetilde f(x,t\,|\,x_0),
 \qquad x\neq x_0,\;\; t>0.
 \label{equation:relgwfw}
\end{equation}
\begin{theorem}
For all $x\in\mathbb R\setminus\{0\}$ and $t>0$ we have
$$
 f(x,t\,|\,0)=\int_0^t f(0,\tau\,|\,0)\,e^{-\nu\,(t-\tau)}\,
 \widetilde g(x,t-\tau\,|\,0)\,{\rm d}\tau.
$$
\end{theorem}
\begin{proof}
We recall that for all $x\in{\mathbb R}$ and $t\geq 0$ the following relation holds:
\begin{equation}
 \widetilde f(x,t\,|\,0)
 =\exp\left\{{2\,(\widehat\lambda-\widehat\mu)\over \sigma^2}\,x\right\}
 \,\widetilde f(0,t\,|\,x).
 \label{equation:idffw}
\end{equation}
Hence, for $x\in{\mathbb R}$ and $t\geq 0$ Eq.\ (\ref{equation:f(x,t)}) becomes
\begin{equation}
 f(x,t\,|\,0)= e^{-\nu t} \,\widetilde f(x,t\,|\,0)
 +\eta \int_0^t q(\tau)\, e^{-\nu\,(t-\tau)}
 \,\exp\left\{{2\,(\widehat\lambda-\widehat\mu)\over \sigma^2}\,x\right\}
 \,\widetilde f(0,t-\tau\,|\,x)\,{\rm d}\tau.
 \label{equation:relffw}
\end{equation}
We recall that for the Wiener process $\widetilde X(t)$ one has:
\begin{equation}
 \widetilde f(0,t\,|\,x)
 =\int_0^t \widetilde g(0,\tau\,|\,x)\,
 \widetilde f(0,t-\tau\,|\,0)\,{\rm d}\tau,
 \qquad x\neq 0,\;\; t>0.
 \label{equation:relfgfw}
\end{equation}
Moreover, making use of (\ref{equation:relgwfw}) and (\ref{equation:idffw}), for 
$x\neq 0$ and $t>0$ we have
\begin{eqnarray}
 \widetilde g(0,t\,|\,x) \!\!\!\!
 &=& \!\!\!\!  
 {|x|\over t}\,\exp\left\{-{2\,(\widehat\lambda-\widehat\mu)\over \sigma^2}\,x\right\}
 \,\widetilde f(x,t\,|\,0)
 \nonumber \\
 &=& \!\!\!\! \exp\left\{-{2\,(\widehat\lambda-\widehat\mu)\over \sigma^2}\,x\right\}
 \,\widetilde g(x,t\,|\,0).
 \label{equation:rlgfg}
\end{eqnarray}
Substituting (\ref{equation:relfgfw}) in the last term of (\ref{equation:relffw}), 
employing Fubini's theorem and making use
again of (\ref{equation:f(x,t)}) we finally obtain:
\begin{eqnarray}
 f(x,t\,|\,0) \!\!\!\!
 &=& \!\!\!\! e^{-\nu t} \,\widetilde f(x,t\,|\,0)
 -e^{-\nu t} \, \exp\left\{{2\,(\widehat\lambda-\widehat\mu)\over \sigma^2}\,x\right\}
 \int_0^t \widetilde g(0,u\,|\,x)\,
 \widetilde f(0,t-u\,|\,0)\,{\rm d}u
 \nonumber \\
 && \!\!\!\! +\exp\left\{{2\,(\widehat\lambda-\widehat\mu)\over \sigma^2}\,x\right\}
 \,\int_0^t \widetilde g(0,u\,|\,x)\, e^{-\nu u} \,f(0,t-u\,|\,0)\,{\rm d}u,
 \qquad x\neq 0.
 \label{eq:pr3.3}
\end{eqnarray}
The first and the second term in the right-hand-side of (\ref{eq:pr3.3}) vanish 
due to Eqs.\ (\ref{equation:idffw}) and (\ref{equation:relfgfw}). 
Hence, the thesis follows making use of (\ref{equation:rlgfg}) in the 
last term of (\ref{eq:pr3.3}). 
\end{proof}
\par
Aiming to analyse the asymptotic behaviour of process $X(t)$ we now 
introduce the steady-state density
\begin{equation}
 W(x)=\lim_{t\to +\infty}f(x,t\,|\,0),
 \qquad x\in{\mathbb R}.
 \label{equation:defW(x)}
\end{equation}
Hereafter we obtain that $W(x)$ is a bilateral asimmetric exponential density.
\begin{theorem}
The steady-state density (\ref{equation:defW(x)}) is given by
\begin{equation}
 W(x)={\eta\,\nu\over\eta+\nu}\,
 {1\over \sqrt{(\widehat\lambda-\widehat\mu)^2 +2\sigma^2 \nu}}\,
 \exp\left\{{\widehat\lambda-\widehat\mu\over\sigma^2}\,x
 -{\sqrt{(\widehat\lambda-\widehat\mu)^2 +2\sigma^2 \nu}\over\sigma^2}\,|x|\right\},
 \qquad x\in{\mathbb R}.
 \label{equation:esprW(x)}
\end{equation}
\end{theorem}
\begin{proof}
Being $W(x)=\displaystyle\lim_{z\downarrow 0}[z\,f^*(x,z\,|\,0)]$,
from (\ref{equation:solphi}) we have (\ref{equation:esprW(x)}).
\end{proof}
\par 
We remark that, due to (\ref{equation:cont}) and (\ref{equation:qstst}), we have 
$$
 \int_{-\infty}^{\infty} W(x)\,{\rm d}x = 1-q = {\eta\over \eta+\nu}.
$$
%   ----------------------------------------------------------------------------
\subsection{Mean and variance}
%   ---------------------------------------------------------------------------- 
Let us now define the mean  
$$
 m_X(t):=E[X(t)\cdot {\bf 1}_{\{X(t)\neq F\}}\,|\,X(0)=0], 
 \qquad t\geq 0.
$$
Due to (\ref{equation:f(x,t)}) and recalling 
that $E[\widetilde X(t)\,|\,\widetilde X(0)=0]=(\widehat\lambda-\widehat\mu)t$, 
for all $t\geq 0$ we have   
$$
 m_X(t)= e^{-\nu t} \,(\widehat\lambda-\widehat\mu)t
 +\eta
 \int_0^t q(\tau)\, e^{-\nu\,(t-\tau)}
 \,(\widehat\lambda-\widehat\mu)(t-\tau)\,{\rm d}\tau.
$$
Use of (\ref{equation:q(t)}) and straightforward calculations thus yield:
\begin{equation}
 m_X(t)
 = \frac{(\widehat\lambda -\widehat\mu)\eta}{(\eta+\nu) \nu}
  \left\{1 - e^{-\nu t}  
 +\frac{\nu^2}{\eta^2} e^{-\nu t}(1 -   e^{-\eta  t}) \right\},
 \qquad t\geq 0.
 \label{eq:medXtfinale}
\end{equation}
The strict correspondence between the means (\ref{eq:medfinale}) and 
(\ref{eq:medXtfinale}) is evident. 
\par
Making use of (\ref{equation:f(x,t)}) one can express the moments of $X(t)$ in terms 
of the moments of the Wiener process. In particular, for any $t\geq 0$ we have: 
\begin{eqnarray}
 && \hspace{-1.2cm} 
  V_{X}(t):=Var[X(t)\cdot {\bf 1}_{\{X(t)\neq F\}}\,|\,X(0)=0] 
  \nonumber \\
 &&   
 = {\sigma^2 \eta\over  (\eta+\nu)\nu}
 \left[ 1-e^{-\nu t} +\frac{\nu^2}{\eta^2} e^{-\nu t}\left(1-e^{-\eta t}\right)\right] 
 \nonumber \\
 &&   +\,{(\widehat\lambda-\widehat\mu)^2\over (\eta+\nu)^2\nu^2\eta^2}
 \Big\{-2\nu^2 e^{-\nu t}\left(1-e^{-\eta t}\right)(\nu^2+\eta\nu+\eta^2)
 +2\nu\eta^3\left(1-e^{-\nu t}\right)
 \nonumber \\
 &&   +\,\eta^4+2\eta\nu(\eta+\nu)(\nu^2-\eta^2)t e^{-\nu t}
 -e^{-2\nu t}\left(\nu^2-\eta^2-\nu^2 e^{-\eta t}\right)^2\Big\}.
 \label{equation:varXt}
\end{eqnarray} 
Notice that for the densities shown (from left to right) in Figure 4 we have 
$m_X(1)=0.865$, $1.305$, $1.834$ and  $V_X(1)=1.283$, $1.469$, $1.198$ (remarking 
that  in this case $V_X(1)$ is not monotonic in $\nu$). The behaviour of mean 
and variance as  $t$ varies is shown in Figure 5. 
\par
Finally, from (\ref{eq:medXtfinale}) and (\ref{equation:varXt}) we immediately obtain the  
asymptotic limits of the mean and variance of the approximating process:
\begin{equation}
 \lim_{t\to +\infty} m_X(t) = \frac{(\widehat\lambda -\widehat\mu)\eta}{(\eta+\nu) \nu},
\qquad 
 \lim_{t\to +\infty} V_X(t)
 ={\sigma^2\eta\over (\eta+\nu)\nu}
 +{(\widehat\lambda-\widehat\mu)^2 \eta\over (\eta+\nu)^2\nu^2}(2\nu+\eta).
\label{eq:limmedvar}
\end{equation}
%
%--------------------------------------------------------------------
\begin{figure}[t]   %%%%%%%%%%% FIGURA  media e varianza di X
%\vspace{20mm}
\begin{center}
\epsfxsize=7.3cm
\epsfbox{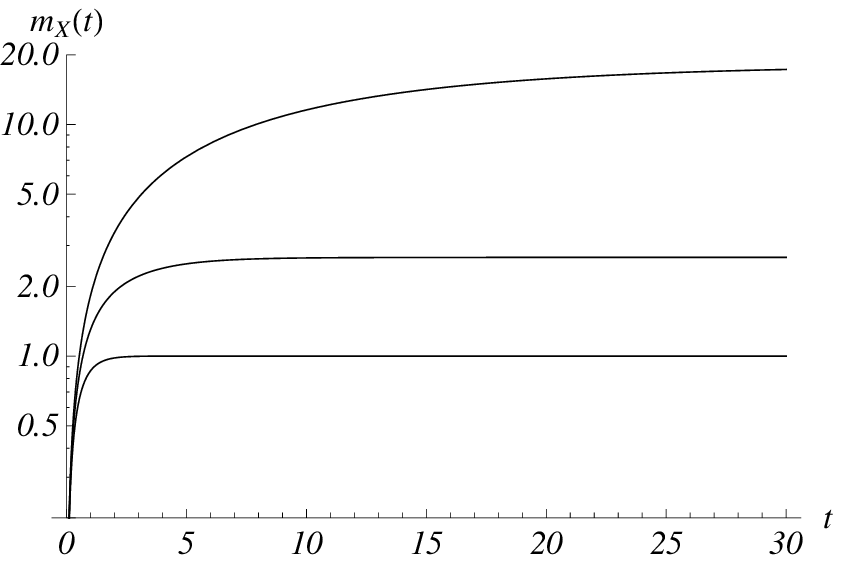} 
$\;\;$
\epsfxsize=7.3cm
\epsfbox{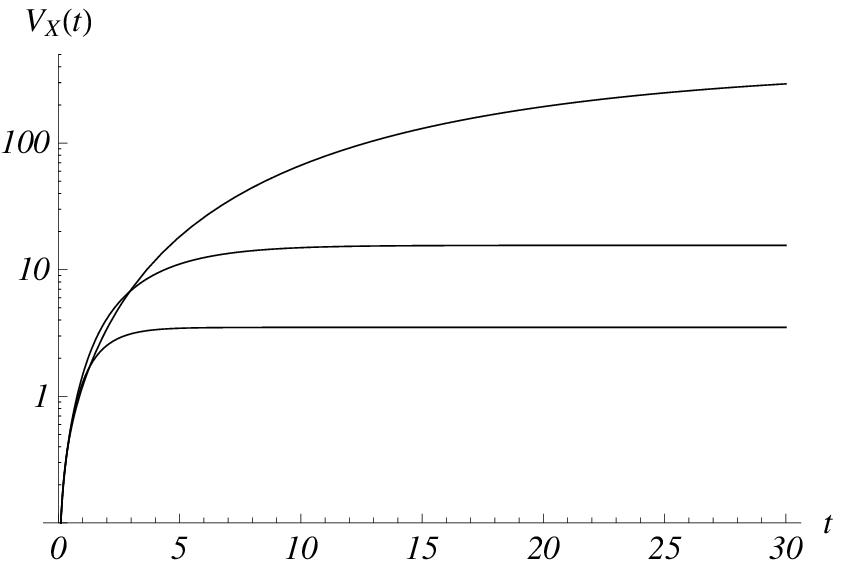}
\vspace{-0.2cm}
\caption{Plots of mean and variance (\ref{eq:medXtfinale}) and (\ref{equation:varXt}) for 
the same values of the parameters as in Figure 4 (from bottom to top).}
\end{center}
\end{figure}
%--------------------------------------------------------------------

%   ----------------------------------------------------------------------------
\section{Remarks on the continuous approximation}
%   ----------------------------------------------------------------------------
\setcounter{equation}{0}
We conclude the paper by giving some comments concerning the goodness 
of the continuous approximation. Recall that, as pointed out in Section 3, 
the approximation is expected to improve as $\varepsilon$ goes to zero 
and as $|n|$ grows larger.
\par
First of all we point out that the means of the discrete process and 
of the continuous approximating process are in strict agreement. Indeed, by 
comparing (\ref{eq:medfinale}) and (\ref{eq:medXtfinale}), under substitutions 
(\ref{equation:parameters}) we immediately have 
$$
 E[N(t)\,|\,N(0)=0]=E\Big[{X(t)\over \varepsilon}\,\Big|\,{X(0)\over \varepsilon}=0\Big]
$$
for all $\varepsilon>0$ and $t\geq 0$. 
\par
In  order to ascertain the approximation between the transient distributions  
of $N(t)$ and $X(t)$, let us now denote by $P_n^*(z;\varepsilon)$, $n\in{\mathbb Z}$,
the Laplace transform of $P_n(t)$ where the substitutions  
(\ref{equation:parameters}) are performed. Hence, making use of
(\ref{equation:pi0}) and (\ref{equation:solphi}) we obtain
$$
 \lim_{\varepsilon\downarrow 0} {P_0^*(z;\varepsilon)\over \varepsilon}
 ={\eta\,\nu\over z\,(z+\eta+\nu)\sqrt{(\widehat\lambda-\widehat\mu)^2+2(z+\nu)\sigma^2}}
 +{1\over \sqrt{(\widehat\lambda-\widehat\mu)^2+2(z+\nu)\sigma^2}}
 =f^*(0,z\,|\,0).
$$
Moreover, by noting that under positions (\ref{equation:parameters})
$$
 \lim_{\varepsilon\downarrow 0}{z+\lambda+\mu+\nu
 -\sqrt{(z+\lambda+\mu+\nu)^2-4\lambda\mu}\over 2\mu}=1,
$$
from (\ref{equation:p*pnz}) for $x\neq 0$ we have
$$
 \lim_{\varepsilon\downarrow 0}{P_n^*(z;\varepsilon)|_{n=x/\varepsilon}\over \varepsilon}
 =\left\{\begin{array}{ll}
  e^{w_2\,x}\,f^*(0,z\,|\,0), & x>0 \\
  e^{w_1\,x}\,f^*(0,z\,|\,0), & x<0 
  \end{array}
  \right\}
  =f^*(x,z\,|\,0),
$$
where $f^*(x,z\,|\,0)$ is given in (\ref{equation:solphi})
and where $w_1$ and $w_2$ are defined in (\ref{equation:wespr}).
Hence, recalling the meaning of $P_n^*(z;\varepsilon)$ and denoting by
$P_n(t;\varepsilon)$ the probabilities (\ref{equation:pn(t)}) when the 
substitutions (\ref{equation:parameters}) are performed, we have:
$$
 \lim_{\varepsilon\downarrow 0}{P_n(t;\varepsilon)|_{n=x/\varepsilon}\over\varepsilon}
 =f(x,t\,|\,0).
$$
This implies that under the approximating procedure the state probabilities of the 
discrete process are close to the density $f(x,t\,|\,0)$ of process $X(t)$. 
In addition, this is confirmed by the data of Table~1 showing that the 
agreement of the two steady-state distributions improves as $\varepsilon$ 
decreases and as $|n|$ increases. Furthermore, we remark that performing the substitutions
\begin{equation}
 \varepsilon_h= \varepsilon\,h, \qquad 
 \widehat\lambda_h= \widehat\lambda\,h, \qquad 
 \widehat\mu_h= \widehat\mu\,h, \qquad 
 \sigma^2_h= \sigma^2\,h, 
 \qquad (h>0)
 \label{eq:sostit}
\end{equation} 
recalling (\ref{equation:parameters}) we have 
$$
 \lambda_h:={\widehat \lambda_h\over \varepsilon_h}+{\sigma^2_h\over 2\varepsilon^2_h}=\lambda,
 \qquad
 \mu_h:={\widehat \mu_h\over \varepsilon_h}+{\sigma^2_h\over 2\varepsilon^2_h}=\mu.
$$
Hence, carrying out substitutions (\ref{eq:sostit}) in the steady-state density 
(\ref{equation:esprW(x)}) we obtain   
$W_h(n  \varepsilon_h) \varepsilon_h=W(n  \varepsilon) \varepsilon$ for any $h>0$. 
This shows that the approximation procedure is effective when the rates 
$\lambda$ and $\mu$ are suitably varying. 
\par
In conclusion we point out that the agreement between the discrete model and the 
approximating jump-diffusion process is also confirmed by the comparison of their 
asymptotic means and variances. Indeed, the correspondence of limits (\ref{eq:limmedvar}) 
with those obtained in Eqs.\ (\ref{eq:medasfinale}) and (\ref{equation:varN}) for the 
discrete model is immediate under positions (\ref{equation:parameters}). 
%
%--------------------- TABELLA UNICA -----------------------------------------------
\begin{table}[t]
\begin{center}
 \small   
\begin{tabular}{|r|ccc|ccc|ccc|}
\hline 
{} & {} & $\varepsilon=0.1$ & {} & {} & $\varepsilon=0.05$  & {} & {} & $\varepsilon=0.01$  & {} \\
\hline 
$n$ & $\pi_n(\varepsilon)/\varepsilon$ & $W(n\varepsilon)$ & $\Delta (n)$ & 
$\pi_n(\varepsilon)/\varepsilon$ & $W(n\varepsilon)$ & $\Delta (n)$ & 
$\pi_n(\varepsilon)/\varepsilon$ & $W(n\varepsilon)$ & $\Delta (n)$  \\
\hline
$-6$ & 0.03621 & 0.03668  & 0.01305  & 0.04073  & 0.04102   & 0.00721 & 0.04480 & 0.04487  & 0.00155  \\
$-5$ & 0.03757 & 0.03807  & 0.01353  & 0.04149  & 0.04180   & 0.00733 & 0.04496 & 0.04503  & 0.00156  \\
$-4$ & 0.03897 & 0.03952  & 0.01401  & 0.04227  & 0.04258   & 0.00745 & 0.04513 & 0.04520  & 0.00156  \\
$-3$ & 0.04044 & 0.04102  & 0.01448  & 0.04306  & 0.04339   & 0.00757 & 0.04530 & 0.04537  & 0.00157  \\
$-2$ & 0.04196 & 0.04258  & 0.01497  & 0.04387  & 0.04420   & 0.00769 & 0.04547 & 0.04554  & 0.00157  \\
$-1$ & 0.04353 & 0.04420  & 0.01544  & 0.04469  & 0.04503   & 0.00781 & 0.04564 & 0.04571  & 0.00158  \\
$0$ & 0.04516  & 0.04588  & 0.01593  & 0.04552   & 0.04588   & 0.00793 & 0.04581  & 0.04588  & 0.00158 \\
$1$ & 0.04260  & 0.04323  & 0.01472  & 0.04420   & 0.04454   & 0.00763 & 0.04554  & 0.04561  & 0.00157 \\
$2$ & 0.04019  & 0.04073  & 0.01351  & 0.04292   & 0.04323   & 0.00732 & 0.04527  & 0.04534  & 0.00156 \\
$3$ & 0.03791  & 0.03838  & 0.01231  & 0.04167   & 0.04196   & 0.00702 & 0.04500  & 0.04507  & 0.00154 \\
$4$ & 0.03576  & 0.03616  & 0.01111  & 0.04046   & 0.04073   & 0.00672 & 0.04473  & 0.04480  & 0.00153 \\
$5$ & 0.03373  & 0.03407  & 0.00990  & 0.03929   & 0.03954   & 0.00641 & 0.04447  & 0.04454  & 0.00152 \\
$6$ & 0.03182  & 0.03210  & 0.00870  & 0.03814   & 0.03838   & 0.00611 & 0.04421  & 0.04427  & 0.00151 \\
\hline
\end{tabular}
\caption{For $\eta=0.25$, $\nu=1$, $\widehat\lambda=1$, $\widehat\mu=2$ and $\sigma=3$   
the quantities $\pi_n(\varepsilon)/\varepsilon$ and $W(n\varepsilon)$ 
are listed together with the relative difference 
$\Delta(n):=[W(n\varepsilon)\varepsilon-\pi_n(\varepsilon)]/\pi_n(\varepsilon)$ 
for various integer values of $n$ and various choices of $\varepsilon$. 
We have  $(\lambda,\mu)=(460,470)$ when $\varepsilon=0.1$, 
$(\lambda,\mu)=(1820,1840)$ when $\varepsilon=0.05$,  and 
$(\lambda,\mu)=(45100,45200)$ when $\varepsilon=0.01$.} 
\end{center}
\end{table}
\small
%
%--------------------------------------------------------------------
\subsection*{\bf Acknowledgments}
%--------------------------------------------------------------------
%
The constructive criticism of two anonymous reviewers is 
gratefully acknowledged. 
\par
This work has been partially performed during a visit of Dr B.\ 
Krishna Kumar at Salerno University. He would like to acknowledge 
the Department of Mathematics and Informatics for kindness and hospitality.
The support by MIUR (PRIN 2008) and by GNCS-INdAM  is 
acknowledged by the other authors.
% 
%--------------------------------------------------------------------

%

\begin{thebibliography}{99}
%--------------------------------------------------------------------
%
%
\bibitem{Aba91}
Abate J, Kijima M, Whitt W (1991) 
Decomposition of the $M/M/1$ transition function. 
Queueing Systems 9:323--336 
%
\bibitem{BaMa89}
Baccelli F, Massey WA (1989)
A sample path analysis of the $M/M/1$ queue. 
J Appl Prob 26:418--422 
%
\bibitem{Br85}
Brockwell PJ (1985) 
The extinction time of a birth, death and 
catastrophe process and of a related diffusion model.
Adv Appl Prob 17:42--52 
%
\bibitem{BDGNR09} 
Buonocore A, Di Crescenzo A, Giorno V, Nobile AG, Ricciardi LM (2009) 
A Markov chain-based model for actomyosin dynamics.  
Scientiae Math Japon 70:159--174  
%
\bibitem{CZ2003} 
Chao X, Zheng Y (2003) 
Transient analysis of immigration-birth-death process with total catastrophes. 
Prob Eng Inform Sci 17:83--106 
%
\bibitem{Co71}
Conolly BW (1971) 
On randomized random walks.
SIAM Review 13:81--99 
%
\bibitem{CoPaSe02}
Conolly BW, Parthasarathy PR, Selvaraju N  (2002)
Doubled-ended queues with impatience. 
Computers Oper Res 29:2053--2072 
%
\bibitem{DiCrGiNoRi03}
Di Crescenzo A, Giorno V, Nobile AG, Ricciardi LM (2003) 
On the $M/M/1$ queue with catastrophes and its continuous approximation. 
Queueing Systems 43:329--347 
%
\bibitem{DiCrGiNoRi08}
Di Crescenzo A, Giorno V, Nobile AG, Ricciardi LM (2008) 
A note on birth-death processes with effective catastrophes. 
Stat Prob Lett 78:2248--2257 
%
\bibitem{DiCrNa04} 
Di Crescenzo A, Nastro A (2004)
On first-passage-time densities for certain symmetric Markov chains.  
Scientiae Math Japon 60:381--390  
%
\bibitem{DiCrNo95}
Di Crescenzo A, Nobile AG (1995)
Diffusion approximation to a queueing system with time-dependent
arrival and service rates. 
Queueing Systems 19:41--62 
%
\bibitem{EcFa2003}
Economou A, Fakinos D (2003) 
A continuous-time Markov chain under the influence of a regulating 
point process and applications in stochastic models with catastrophes. 
Europ J Operat Res 149:625--640    
%
\bibitem{EcFa2008}
Economou A, Fakinos D (2008) 
Alternative approaches for the transient analysis 
of Markov chains with catastrophes. 
J Statist Theory Pract 2:183--197   
%
\bibitem{GaSw2006}
Gani J, Swift RJ (2006) 
A simple approach to birth processes with random catastrophes. 
J Combin Inform System Sci 31:1--7  
%
\bibitem{GaSw2006}
Gani J, Swift RJ (2007) 
Death and birth-death and immigration processes with catastrophes. 
J Statist Theory Pract 1:39--48   
%
\bibitem{GiNoRi86}
Giorno V, Nobile AG, Ricciardi LM (1986) 
On some diffusion approximations to queueing systems. 
Adv Appl Prob 18:991--1014 
%
\bibitem{GiNoRi87}
Giorno V, Nobile AG, Ricciardi LM (1987) 
On some time non-homogeneous diffusion approximations to queueing systems. 
Adv Appl Prob  19:974--994 
%
\bibitem{JMB2007}
Jain JL, Mohanty SG, Boehm W (2007) 
A Course on Queueing Models. 
Chapman and Hall/CRC, Boca-Raton   
%
\bibitem{KeWi1990}
Kella O, Whitt W (1990)
Diffusion approximation for queues with server vacations.
Adv Appl Prob  22:706--729  
%
\bibitem{Ki2004}
Kimura T (2004)
Diffusion models for computer/communication systems. 
Econ J  Hokkaido Univ 33:37--52  
%
\bibitem{KrMo2007}
Krinik A, Mortensen C (2007) 
Transient probability functions of finite birth-death processes with catastrophes. 
J Statist Plann Infer 137:1530--1543  
%
\bibitem{Kretal2005}
Krinik A, Rubino G, Marcus D, Swift R, Kasfy H, Lam H (2005) 
Dual processes to solve single server systems. 
J Statist Plann Infer  135:121--147   
%
\bibitem{KuAr2000}
Krishna Kumar B, Arivudainambi D (2000)
Transient solution of an $M/M/1$ queue with catastrophes.
Comput Math Appl 40:1233--1240 
%
\bibitem{KuMa2003}
Krishna Kumar B, Pavai Madheswari S (2003)
Transient solution of an $M/M/2$ queue with catastrophes.
Math Scientist 28:129--136 
%
\bibitem{KuKrMaBa2007}
Krishna Kumar B, Krishnamoorthy A, Pavai Madheswari S, Sadiq Basha S (2007) 
Transient analysis of a single server queue with
catastrophes, failures and repairs. 
Queueing Syst 56:133--141 
%
\bibitem{PaLe2004}
Parthasarathy PR, Lenin RB (2004) 
Birth and death processes (BDP) models with applications---queueing  
communication systems, chemical models. Biological models: 
the state-of-the-art with a time-dependent perspective. 
Amer J Math Management Sci 24:1--212   
%
\bibitem{St2006}
Stirzaker D (2006) 
Processes with catastrophes. 
Math Scientist 31:107--118   
%
\bibitem{St2006}
Stirzaker D (2007) 
Processes with random regulation. 
Prob Eng Inform Sci 21:1--17   
%
\bibitem{St2006}
Swift RJ (2000) 
A simple immigration-catastrophe process. 
Math Scientist 25:32--36   
%
\bibitem{Sw01}
Swift RJ (2001) 
Transient probabilities for a simple birth-death-immigration process under 
the influence of total catastrophes.
Int J Math Math Sci 25:689--692
%
\bibitem{Sw04}
Switkes J (2004) 
An unbiased random walk with catastrophe. 
Math Scientist 29:115--121 
%
\end{thebibliography}
\end{document}